\documentclass[a4,12pt,multicol,makeidx]{article}
\def\origin{
\hbox{}\vskip-\baselineskip\vskip-\topskip%
  \vbox to 0pt{\vskip-1in%
    \hbox to 0pt{\hskip-1in%
      \hbox to 0pt{\vrule width 1cm height .4pt depth 0mm\hss}%
      \vbox to 0pt{\hrule width .4pt height 0pt depth 1cm\vss}%
    \hss}%
  \vss}
  \vskip-\baselineskip
  \vbox to 0pt{\vskip-1in\vskip3cm%
    \hbox to 0pt{\hskip-1in\hskip3cm%
      \hbox to 0pt{\hss\vrule width 2cm height .4pt depth 0mm\hss}%
      \vbox to 0pt{\vss\hrule width .4pt height 1cm depth 1cm\vss}%
    \hss}%
  \vss}%
  \vskip5mm\hskip10mm (3cm,3cm)
}%

\def\a{\alpha} \def\d{\delta}    \def\p{\partial}  \def\e{\varepsilon} 
\def\n{\nabla}   
 
\def\leq{\underline{<}} 
 
\def\ox{\overline{x}} \def\oy{\overline{y}}  
\def\hx{\hat{x}} \def\hy{\hat{y}} 
\def\la{\langle} \def\ra{\rangle}
\def\ox{\overline{x}}  
\def\oy{\overline{y}}

\newenvironment{theorem}{%
\par \bigskip \it}{%
\bigskip \par}

\newpage
\makeindex
\title{Corrigendum for the comparison theorems in \\
"A new definition of viscosity solutions 
for a class of second-order degenerate elliptic 
integro-differential equations".}
\author{Mariko Arisawa\\ GSIS, Tohoku University
\\ Aramaki 09, Aoba-ku, Sendai 980-8579, JAPAN\\
E-mail: arisawa@math.is.tohoku.ac.jp
}
\date{}
\begin{document}
\maketitle
\bigskip

	In this note, we shall present the correction of the proofs of 
the comparison results in the paper [1]. In order to show clearly the correct way of the demonstration, we shall simplify the problem to the following.\\
$$\hbox{(Problem (I))}:\quad
	F(x,u,\nabla u,\nabla^2 u) -\int_{{\bf R^N}}
	u(x+z)-u(x)\qquad \qquad\qquad\qquad
$$
\begin{equation}\label{problem}
	\qquad -{\bf 1}_{|z|\leq 1}\la z,\n u(x) \ra q(dz)=0 \qquad \hbox{in} \quad \Omega,
\end{equation}

$$\hbox{(Problem (II))}:\quad
	F(x,u,\nabla u,\nabla^2 u) -\int_{\{z\in {\bf R^N}|x+z\in \overline{\Omega}\}}
	u(x+z)-u(x)\qquad \qquad
$$
\begin{equation}\label{problem2}
	\qquad -{\bf 1}_{|z|\leq 1}\la z,\n u(x) \ra q(dz)=0 \qquad \hbox{in} \quad \Omega,
\end{equation}
where $\Omega\subset {\bf R^N}$ is open, and  $q(dz)$ is a positive Radon measure such that 
$\int_{|z|\leq 1} |z|^2 q(dz)+\int_{|z| > 1} 1 q(dz) <\infty$. 
 Although in [1] only (II) 
was studied,  in order to avoid the non-essential  technical complexity, here, let us give the explanation mainly for (I).  For (I), we consider the Dirichlet B.C.:
\begin{equation}\label{dirichletI}
	u(x)=g(x) \quad \forall x\in \Omega^c,
\end{equation}
where $g$ is a given continuous function in $\Omega^c$. 
 For (II), we assume that $\Omega$ is a precompact convex open subset in ${\bf R^N}$ with $C^1$ boundary satisfying the uniform exterior sphere condition, and consider either the Dirichlet B.C.:
 \begin{equation}\label{dirichletII}
	u(x)=h(x) \quad \forall x\in \p\Omega,
\end{equation}
where $h$ is a given continuous function on $\p\Omega$, or the Neumann B.C.:
\begin{equation}\label{neumann}
	\la \n u(x), {\bf n}(x)\ra=0 \quad \forall x\in \p\Omega,
\end{equation}
where ${\bf n}(x)\in {\bf R^N}$ the outward unit normal vector field defined on $\p\Omega$. The above problems are studied in the framework of the viscosity solutions introduced in [1]. Under all the assumptions in [1], for (I) the following comparison result holds, and for (II), although the proofs therein are incomplete, the comparison results stated in [1] hold, and we shall show in a future article.\\
\begin{theorem}{\bf Theorem 1.1\quad (Problem I with Dirichlet B.C.)}  Assume that $\Omega$ is bounded, and the conditions for $F$ in [1] hold. Let $u\in USC({\bf R^N})$ and $v\in LSC(\bf R^N)$ be respectively a viscosity subsolution and a supersolution of (\ref{problem}) in $\Omega$, which satisfy $u\leq v$ on $\Omega^c$. Then, 
$u\leq v$ in $\Omega$. 
\end{theorem}

	To prove Theorem 1.1, we approximate the solutions $u$ and $v$ by the supconvolution: $u^r(x)=\sup_{y\in {\bf R^N}}\{ 
	u(y)-\frac{1}{2r^2}|x-y|^2\}$ and the infconvolution: 
	$v_r(x)=\inf_{y\in {\bf R^N}}\{ 
	v(y)+\frac{1}{2r^2}|x-y|^2
	\}$ $(x\in {\bf R^N})$, where $r>0$. \\
	
\begin{theorem} {\bf Lemma 1.2\quad (Approximation for Problem (I))} Let $u$ and $v$ be respectively a viscosity subsolution and a  supersolution of (\ref{problem}). For any $\nu>0$ there exists $r>0$ such that $u^r$ and $v_r$ are respectively a subsolution and a supersolution of the following problems. 
\begin{equation}\label{subproblem}
	F(x,u,\nabla u,\nabla^2 u) -\int_{{\bf R^N}}
	u(x+z)-u(x)-{\bf 1}_{|z|\leq 1}\la z,\n u(x) \ra q(dz)\}\leq \nu,\quad
\end{equation}
\begin{equation}\label{superproblem}
	F(x,v,\nabla v,\nabla^2 v) -\int_{{\bf R^N}}
	v(x+z)-v(x)-{\bf 1}_{|z|\leq 1}\la z,\n v(x) \ra q(dz)\}\geq -\nu,
\end{equation}
in $\Omega_r=\{x\in {\Omega}|\quad dist(x,\p\Omega)> \sqrt{2M}r\}$, where $M=\max\{\sup_{\overline{\Omega}}|u|,\sup_{\overline{\Omega}}|v|\}$. 
\end{theorem}

Remark that $u^r$ is semiconvex, $v_r$ is semiconcave, and both are Lipschitz continuous in ${\bf R^N}$. We then deduce from the Jensen's maximum principle and the Alexandrov's theorem (deep results in the convex analysis, see [2] and [3]), the following lemma, the last claim of which is quite important in the limit procedure in the nonlocal term. 
\begin{theorem}{\bf Lemma 1.3}
	Let $U$ be semiconvex and $V$ be semiconcave in $\Omega$. 
	For 
$\phi(x,y)=\a |x-y|^2$ ($\a>0$) consider 
$\Phi(x,y)=U(x)-V(y)-\phi(x,y)$, 
and assume that $(\ox,\oy)$ is an interior maximum of $\Phi$ 
in $\overline{\Omega}\times \overline{\Omega}$. Assume also that there is an open precompact subset $O$ of $\Omega\times \Omega$ containing $(\ox,\oy)$,  and that 
$\mu$
$=\sup_{O} \Phi(x,y) -\sup_{\p O} \Phi(x,y) >0$.
Then, the following holds.\\
(i) There exists a sequence of points $(x_m,y_m)\in O$ ($m\in {\bf N}$) 
such that $\lim_{m\to \infty} (x_m,y_m)=(\ox,\oy)$, and 
$(p_m,X_m)\in J^{2,+}_{\Omega}U(x_m)$,  $(p_m',Y_m)\in J^{2,-}_{\Omega}V(y_m)$
 such that 
$\lim_{m\to \infty}p_m$$=\lim_{m\to \infty}p'_m$$=2\a(x_m-y_m)=p$, and 
$X_m\leq Y_m\quad \forall m$.\\
(ii) For $P_m=(p_m-p,-({p'}_m-p))$, $\Phi_m(x,y)=\Phi(x,y)-\la P_m,(x,y)\ra$
 takes a maximum at $(x_m,y_m)$ in O. \\
(iii) The following holds for any $z\in {\bf R^N}$ such that
 $(x_m+z,y_m+z)\in O$.
\begin{equation}\label{important}
	U(x_m+z)-U(x_m)-\la p_m,z\ra
	\leq V(y_m+z)-V(y_m)-\la p'_m,z\ra. 
\end{equation}
\end{theorem}

	By admitting these lemmas here, let us show how Theorem 1.1 
is proved. \\
$Proof\quad of\quad Theorem\quad 1.1.$ 
We use the argument by contradiction, and assume that $\max_{\overline{\Omega}}(u-v)$$=(u-v)(x_0)=M_0>0$ for $x_0\in \Omega$. Then, we approximate $u$ by $u^r$ (supconvolution) and $v$ by $v_r$ (infconvolution), which are a subsolution  and a supersolution of (\ref{subproblem}) and (\ref{superproblem}), respectively. Clearly, $\max_{\overline{\Omega}}(u^r-v_r)$$\geq M_0>0$. Let $\ox\in \Omega$ be the maximizer of $u^r-v_r$. In the following, we abbreviate the index and  write $u=u^r$, $v=v_r$ without any confusion. As in the PDE theory, consider $\Phi(x,y)=u(x)-v(y)-\a|x-y|^2$, and let $(\hx,\hy)$  be the maximizer of $\Phi$. Then, from Lemma 1.3 there exists $(x_m,y_m)\in \Omega$ ($m\in {\bf N}$) 
such that $\lim_{m\to \infty} (x_m,y_m)=(\hx,\hy)$, and we can take 
$(\e_m,\d_m)$  a pair of positive numbers such that
$u(x_m+z)\leq u(x_m)+\la p_m,z\ra+\frac{1}{2}\la X_m z,z\ra+\d_m|z|^2$, 
$v(y_m+z)\geq v(y_m)+\la p'_m,z\ra+\frac{1}{2}\la Y_m z,z\ra-\d_m|z|^2$, 
for $\forall |z|\leq \e_m$. 
	From the definition of the viscosity solutions, we have 
$$
	F(x_m,u(x_m),p_m,X_m) -\int_{|z|\leq \e_m}
	\frac{1}{2}\la (X_m+2\d_m I)z,z \ra dq(z)
$$
$$
	-\int_{|z|\geq \e_m}
	u(x_m+z)-u(x_m)
	-{\bf 1}_{|z|\leq 1}\la z,p_m \ra q(dz)\leq \nu,
$$
$$
	F(y_m,v(y_m),p'_m,Y_m) -\int_{|z|\leq \e_m}
	\frac{1}{2}\la (Y_m-2\d_m I)z,z \ra dq(z)
$$
$$
	-\int_{|z|\geq \e_m}
	v(y_m+z)-v(y_m)
	-{\bf 1}_{|z|\leq 1}\la z,p'_m \ra q(dz)\geq -\nu.
$$	
	By taking the difference of the above two inequalities,  by using (\ref{important}), and by passing $m\to \infty$ (thanking to (\ref{important}), it is now available), we can obtain the desired contradiction. The claim 
$u\leq v$ is proved. \\

{\bf Remark 1.1.} To prove the comparison results for (II) (in [1]), we do 
 the approximation by the supconvolution: $u^r(x)=\sup_{y\in {\overline{\Omega}}}\{ {u}(y)-\frac{1}{2r^2}|x-y|^2\}$, and the infconvolution: $v_r(x)=\inf_{y\in \overline{\Omega}}\{ {v}(y)+\frac{1}{2r^2}|x-y|^2\}$ as in Lemma 1.2. Because of the restriction of the domain of the integral of the nonlocal term and the Neumann B.C.,  a slight technical complexity is added. The approximating problem for (\ref{problem2})-(\ref{neumann}) in $\overline{\Omega}$ is as follows. 
$$
	\min[F(x,u(x),\nabla u(x),\nabla^2 u(x)) + 
	\min_{y\in \overline{\Omega},|x-y|\leq \sqrt{2M}r }\{
	-\int_{\{z\in \bf R^N|y+z\in \overline{\Omega}\}}
	u(x+z)-u(x)
$$
$$
	-{\bf 1}_{|z|\leq 1}\la z,\n u(x) \ra q(dz),\quad \min_{y\in \p{\Omega},|x-y|\leq \sqrt{2M}r}\{ \la {\bf n}(y),\n u(x) \ra+\rho\} ]
	\leq \nu,
$$
$$
	\max[F(x,v(x),\nabla v(x),\nabla^2 v(x)) + 
	\max_{y\in \overline{\Omega},|x-y|\leq \sqrt{2M}r}\{
	-\int_{\{z\in \bf R^N|y+z\in \overline{\Omega}\}}
	v(x+z)-v(x)
$$
$$
	-{\bf 1}_{|z|\leq 1}\la z,\n v(x) \ra q(dz),\quad \max_{y\in \p{\Omega},|x-y|\leq \sqrt{2M}r}\{ \la {\bf n}(y),\n v(x) \ra-\rho\}]
	\geq -\nu.
$$
We deduce the comparison result from this approximation and Lemma 1.3, by using  the similar argument as in the proof of Theorem 1.1.\\


\end{document}